\ifx\fltex\undefined
\documentclass[11pt,a4paper,english]{article}
\usepackage{amsmath,amssymb}
\usepackage{babel}
\usepackage{eucal}
\fi 

\DeclareMathVersion{scriptstyle}

\usepackage{isolatin1,theorem}
\theoremstyle{change}
{\theorembodyfont{\slshape} 
  \newtheorem{thm}{Theorem.}[section]
  \newtheorem{lemma}[thm]{Lemma.}
  \newtheorem{prop}[thm]{Proposition.}
  \newtheorem{cor}[thm]{Corollary.}
}
{\theorembodyfont{\rmfamily} 

}

\numberwithin{equation}{section} 

\newcommand{\proof}[1][Proof. ]{{\it#1}}

\def\endproof{{\nobreak\qquad{\scriptstyle \blacksquare}}}

\def\C{{\mathbb C}}

\def\N{{\mathbb N}}
\def\R{{\mathbb R}}

\def\CH{{\mathcal H}}

\def\CJ{{\mathcal J}}

\def\CM{{\mathcal M}}

\def\CR{{\mathcal R}}

\def\unit{{\bf 1}}

\def\d{{\rm d}}
\def\eps{\varepsilon}
\def\<{{\langle}}
\def\>{{\rangle}}

\def\Prob{{\text{\rm Prob}}}

\def\cc{^*}

\begin{document}

\title{{\bf{\LARGE{Semicircularity, Gaussianity and\\ Monotonicity of Entropy}}}}
\author{Hanne~Schultz \footnote{As a student of the
      PhD-school OP-ALG-TOP-GEO the author is partially supported by the Danish Research Training Council.} \footnote{Partially
      supported by The Danish National Research Foundation.}}
\date{}       
\maketitle

\begin{center}
{\it 
Para el Grupo
}
\end{center}

\begin{abstract}
\noindent S.~Artstein, K.~Ball, F.~Barthe, and A.~Naor have shown (cf. \cite{ABBN}) that if $(X_j)_{j=1}^\infty$ are i.i.d. random variables, then the entropy of ${\textstyle \frac{X_1+\cdots+X_{n}}{\sqrt{n}}}$, $H\Big({\textstyle \frac{X_1+\cdots+X_{n}}{\sqrt{n}}}\Big)$, increases as $n$ increases. The free analogue 
was recently proven by D.~Shlyakhtenko in \cite{S}. That is, if $(x_j)_{j=1}^\infty$ are freely independent, identically distributed, self-adjoint elements in a noncommutative probability space, then the free entropy of ${\textstyle \frac{x_1+\cdots+x_{n}}{\sqrt{n}}}$, $\chi\Big({\textstyle \frac{x_1+\cdots+x_{n}}{\sqrt{n}}}\Big)$, increases as $n$ increases. In this paper we prove that if $H(X_1)>-\infty$ ($\chi(x_1)>-\infty$, resp.), and if the entropy (the free entropy, resp.) is {\it not} a strictly increasing function of $n$, then $X_1$ ($x_1$, resp.) must be Gaussian (semicircular, resp.). 
\end{abstract}

\section{Introduction.}

Shannon's entropy of a (classical) random variable $X$ with Lebesgue absolutely continuous distribution $\d\mu_X(x)=\rho(x)\d x$, is given by
\begin{equation}
H(X)=-\int_\R \rho(x)\log \rho(x)\d x,
\end{equation}
whenever the integral exists. If the integral does not exist, or if the distribution of $X$ is not Lebesgue absolutely continuous, then $H(X)=-\infty$. 

The entropy can also be written in terms of score functions and of Fisher information. Take a standard Gaussian random variable $G$ such that $X$ and $G$ are independent. Let 
\[
X^{(t)}= X+\sqrt t G, \qquad t\geq 0,
\]
and let $j(X^{(t)})=\big(\frac{\partial}{\partial x}\big)\cc (\unit)\in L^2(\mu_{X^{(t)}})$ denote the score function of $X^{(t)}$ (cf. \cite[Section~3]{S}). Then
\begin{equation}\label{entropy, score}
H(X)= \frac12 \int_0^\infty \Big[{\frac{1}{1+t}}-\|j(X^{(t)})\|_2^2\Big]\d t +\frac12 \log (2\pi e). 
\end{equation}
The quantity $\|j(X^{(t)})\|_2^2$ is called the Fisher information of $X^{(t)}$ and is denoted by $F(X^{(t)})$. Among all random variables with a given variance, the Gaussians are the (unique) ones with the smallest Fisher information and the largest entropy.  

A.~J.~Stam (cf. \cite{St}) was the first to rigorously show that if $X_1$ and $X_2$ are independent random variables of the same variance, with $H(X_1), H(X_2)>-\infty$, then for all $t\in [0,1]$,
\[
H(\sqrt t X_1 + \sqrt{1-t}X_2) \geq tH(X_1)+ (1-t)H(X_2),
\]
with equality iff $X_1$ and $X_2$ are Gaussian. It follows that if $(X_j)_{j=1}^\infty$ is a sequence of i.i.d. random variables with finite entropy, then 
\[
n\mapsto H\Big({\textstyle \frac{X_1+\cdots+ X_{2^n}}{2^{\frac n2}}}\Big)
\]
is an increasing function of $n$, and if it is not {\it strictly}
increasing, then $X_1$ is necessarily Gaussian.

Knowing about Stam's result, it seems natural to ask whether the
map
\[
n\mapsto H\Big({\textstyle \frac{X_1+\cdots+ X_{n}}{\sqrt n}}\Big)
\]
is monotonically increasing as well, or even simpler: Is
$H\Big(\frac{X_1+X_2+X_3}{\sqrt 3}\Big)\geq H\Big(\frac{X_1+X_2}{\sqrt
  2}\Big)$? Surprisingly enough, it took more than 40 years for
someone to answer these questions. Both questions were answered in the
affirmative in \cite{ABBN} in 2004.

In this paper we extend Stam's
result by showing that if $H(X_1)>-\infty$ and if for some $n\in\N$, 
\[
H\Big({\textstyle \frac{X_1+\cdots+ X_{n+1}}{\sqrt {n+1}}}\Big)= H\Big({\textstyle \frac{X_1+\cdots+ X_{n}}{\sqrt n}}\Big),
\]
then $X_1$ is necessarily Gaussian (Theorem~\ref{classical}). 


\vspace{.2cm}

\noindent Free entropy, which is the proper free
analogue of Shannon's entropy, was defined by Voiculescu in
\cite{V1}. If $x$ is a self-adjoint element in a finite von Neumann
algebra $\CM$ with faithful normal tracial state $\tau$ and if
$\mu_x\in \Prob(\R)$ denotes the distribution of $x$ with respect to $\tau$, then the free entropy of $x$, $\chi(x)\in [-\infty,\infty[$, is given by
\[
\chi(x) = \int\int \log|s-t| \d\mu_x(s) \d\mu_x(t) + \frac34 + \frac 12 \log(2\pi).
\]
Exactly as in the classical case, $\chi(x)$ may be written in terms of the free analogue of the score function (the conjugate variable) and the free Fisher information. That is, if $s$ is a (0,1)-semicircular element which is freely independent of $x$ and if we let 
\[
x^{(t)}= x+\sqrt t s, \qquad t\geq 0,
\]
then
\begin{equation}\label{free entropy, score}
\chi(x)= \frac12 \int_0^\infty\Big[{\frac{1}{1+t}}-\Phi(x^{(t)})\Big]\d t +\frac12 \log (2\pi e), 
\end{equation} 
where $\Phi(x^{(t)})$ is the free Fisher information of $x^{(t)}$. In \cite{V2} Voiculescu defines for a (non-scalar) self-adjoint variable $y$ in $(\CM,\tau)$ a derivation $\partial_y : \C[y]\rightarrow \C[y]\otimes\C[y]$ by
\[
\partial_y (\unit) =0 \qquad {\rm and} \qquad \partial_y (y) = \unit\otimes\unit.
\]
Then the conjugate variable of $y$, if it exists, is the unique vector $\CJ(y)\in L^2(W\cc(y))$ satisfying that for all $k\in\N$,
\begin{equation}\label{conjugate variable}
\<\CJ(y),y^k\> = \<\unit\otimes\unit, \partial_y (y^k)\>.
\end{equation}
That is, $\CJ(y)=(\partial_y)\cc (\unit\otimes\unit)$. The conjugate variable is the free analogue of the score function, and the free Fisher information of $y$ is exactly $\|\CJ(y)\|_2^2$, so that
\begin{equation}\label{free entropy, score II}
\chi(x)= \frac12 \int_0^\infty\Big[{{\frac{1}{1+t}}}-\|\CJ(x^{(t)})\|_2^2\Big]\d t +\frac12 \log (2\pi e). 
\end{equation} 

Note that if $\CJ(y)=y$, then the moments of $y$ are determined by
\eqref{conjugate variable}, and it is not hard to see that $y$ is necessarily (0,1)-semicircular. 

In \cite{S} D.~Shlyakhtenko showed that if $(x_j)_{j=1}^\infty$ are freely independent, identically distributed self-adjoint elements in $(\CM,\tau)$, then the map
\[
n\mapsto \chi\Big({\textstyle \frac{x_1+\cdots+x_n}{\sqrt n}}\Big)
\]
is monotonically increasing in $n$. In fact, the method used in
\cite{S} applies to the classical case as well. In this paper we will dig into the proof of the inequality
\begin{equation}\label{intro1}
\chi\Big({\textstyle \frac{x_1+\cdots+x_{n+1}}{\sqrt {n+1}}}\Big)
\geq \chi\Big({\textstyle \frac{x_1+\cdots+x_n}{\sqrt n}}\Big)
\end{equation}
and find out what it means for all of the estimates obtained in the course of the proof to be equalities. We conclude that if $\chi(x_1)>-\infty$ and if \eqref{intro1} is an equality for some $n$, then $x_1$ is necessarily semicircular. With a few modifications, our method applies to the classical case as well. 

\vspace{.2cm}

\noindent {\it Acknowledgements.} I would like to thank Dimitri~Shlyakhtenko for introducing me to the problem on monotonicity of entropy and for fruitful discussions about the subject. In addition, I would like to thank the UCLA Math Department, especially the functional analysis group, for their kind hospitality during my stay at UCLA.

\section{The Free Case.}

\noindent Recall that the $(0,1)$-{\it semicircle law} is the Lebesgue
absolutely continuous probability measure on $\R$ with density
\[
\d\sigma_{0,1}(t)= \frac{1}{2\pi} \sqrt{4-t^2} \; 1_{[-2,2]}(t)\;\d t.
\]
More generally, for $\mu,\gamma\in\R$ with $\gamma >0$, the
$(\mu,\gamma)$-{\it semicircle law}  is the Lebesgue
absolutely continuous probability measure on $\R$ with density
\[
\d\sigma_{\mu,\gamma}(t)= \frac{1}{2\pi\gamma}
\sqrt{4\gamma-(t-\mu)^2} \; 1_{[\mu-2\sqrt\gamma ,\mu + 2\sqrt\gamma]}(t)\;\d t.
\]
The parameters $\mu$ and $\gamma$ refer to the first moment and the
variance of $\sigma_{\mu,\gamma}$, respectively. 

Throughout this section, $\CM$ denotes a finite von Neumann algebra with faithful, normal, tracial state $\tau$. We are going to prove:

\begin{thm}\label{thm2} Let $n\in\N$ and let $x_1,\ldots,x_{n+1}$ be freely independent, identically distributed self-adjoint elements in $(\CM,\tau)$. Then 
\begin{equation}\label{eq2}
\chi\Big(\frac{x_1+\cdots+x_{n+1}}{\sqrt{n+1}}\Big) \geq \chi\Big(\frac{x_1+\cdots+x_n}{\sqrt n}\Big).
\end{equation}
Moreover, if $\chi(x_1)>-\infty$, then equality holds in \eqref{eq2} iff $x_1$ is semicircular. 
\end{thm}

\noindent Monotonicity of free entropy was already proven in \cite{S}. Likewise, most of the results stated in this section consist of two parts: An inequality which was proven in \cite{S} or in \cite{ABBN} and a second part which was proven by us. 

\begin{prop}\label{thm1} Let $n\in\N$ and let $x_1,\ldots,x_{n+1}$ be freely independent self-adjoint elements in $(\CM,\tau)$ with $\tau(x_j)=0$ and $\|x_j\|_2=\|x_1\|_2$, $1\leq j\leq n+1$. Let $a_1,\ldots,a_{n+1}\in\R$ with $\sum_j a_j^2 =1$, and let $b_1,\ldots,b_{n+1}\in\R$ such that $\sum_j b_j\sqrt{1-a_j^2}=1$. Then
\begin{equation}\label{eq1}
\Phi\Big(\sum_{j=1}^{n+1}a_jx_j\Big) \leq n \sum_{j=1}^{n+1} b_j^2 \Phi\Big({\textstyle \frac{1}{\sqrt{1-a_j^2}}}\sum_{i\neq j}a_ix_i\Big).
\end{equation}
Moreover, if $\Phi(\sum_{i\neq j}a_ix_i)$ is finite for all $j$, then equality in \eqref{eq1} implies that 
\begin{equation}\label{eq3}
\CJ\Big({\textstyle \frac{1}{\|x_1\|_2}}\sum_{j=1}^{n+1}a_jx_j\Big) = {\textstyle \frac{1}{\|x_1\|_2}}\sum_{j=1}^{n+1}a_jx_j,
\end{equation}
so that $\sum_{j=1}^{n+1}a_jx_j$ is $(0,\|x_1\|_2^2)$-semicircular.
\end{prop}

\begin{lemma}\label{commuting projections} Let $P_1,\ldots,P_m$ be commuting projections on a Hilbert space $\CH$. If $\xi_1,\ldots,\xi_m\in\CH$ satisfy that for all $1\leq i\leq m$,
\[
P_1P_2\cdots P_m \xi_i=0,
\]
then 
\begin{equation}\label{eq6}
\|P_1\xi_1+\ldots + P_m\xi_m\|^2 \leq (m-1)\sum_{i=1}^m \|\xi_i\|^2.
\end{equation}
Moreover, if equality holds in \eqref{eq6}, then $\xi_i\in \bigoplus_{j\neq i} \CH_j$,
where
\[
\CH_j := \{\xi\in\CH\,|\, P_k\xi = \xi, \;k\neq j, \; P_j\xi =0\}=
\Big(\bigcap_{k\neq j}P_k(\CH)\Big)\bigcap P_j^\bot (\CH).
\]
\end{lemma}

\noindent \proof The inequality \eqref{eq6} is the content of \cite[Lemma~5]{ABBN}. The starting point of their proof is to write each $\xi_i$ as an orthogonal sum,
\[
\xi_i=\sum_{\eps\in\{0,1\}^m\setminus (1,1,\ldots ,1)}\xi_\eps^i,
\] 
where for $\eps\in\{0,1\}^m\setminus (1,1,\ldots ,1)$,
\[
\xi_\eps^i\in \CH_\eps := \{\xi\in\CH \,|\, P_j\xi = \eps_j\xi, \; 1\leq j\leq m\}.
\]
Then 
\[
P_1\xi_1+\ldots + P_m\xi_m = \sum_{\eps\in\{0,1\}^m\setminus (1,1,\ldots ,1)}\sum_{ \eps_i=1}P_i\xi_\eps^i,
\]
and
\[
\|P_1\xi_1+\ldots + P_m\xi_m\|^2 = \sum_{\eps\in\{0,1\}^m\setminus (1,1,\ldots ,1)} \Big\|\sum_{\eps_i=1}P_i\xi_\eps^i\Big\|^2.
\]
For fixed $\eps\neq (1,1,\ldots, 1)$ there can be at most $m-1$ $i$'s for which $\eps_i=1$. Thus, by the Cauchy-Schwarz inequality,
\begin{equation}\label{eq12}
\Big\|\sum_{\eps_i=1}P_i\xi_\eps^i\Big\|^2 \leq \Big(\sum_{\eps_i=1}\|P_i\xi_\eps^i\|\Big)^2\leq (m-1) \sum_{\eps_i=1}\|P_i\xi_\eps^i\|^2,
\end{equation}
with the second inequality being an equality iff the vector $(\|P_i\xi_\eps^i\|)_{\eps_i=1}$ ($=(\|\xi_\eps^i\|)_{\eps_i=1}$) has $m-1$ coordinates and is parallel to the vector $v=(1,1,\ldots,1)\in\R^{m-1}$. In particular, if the second inequality in \eqref{eq12} is an equality for some $\eps\in\{0,1\}^m$ with more than one coordinate which is zero, then $(\|P_i\xi_\eps^i\|)_{i=1}^m$ must consist of zeros only. It follows now that
\begin{eqnarray} 
\|P_1\xi_1+\ldots + P_m\xi_m\|^2 &\leq & (m-1)\sum_{\eps\in\{0,1\}^m\setminus (1,1,\ldots ,1)} \sum_{\eps_i=1}\|P_i\xi_\eps^i\|^2 \label{ret1}\\
&=& (m-1)\sum_{\eps\in\{0,1\}^m\setminus (1,1,\ldots ,1)}\sum_{i=1}^m\|P_i\xi_\eps^i\|^2 \label{eq22}\\
&\leq &(m-1) \sum_{i=1}^m \sum_{\eps\in\{0,1\}^m\setminus (1,1,\ldots ,1)} \|\xi_\eps^i\|^2\label{eq21}\label{ret2}\\
&= &(m-1) \sum_{i=1}^m \|\xi_i\|^2.
\end{eqnarray}
Moreover, equality in \eqref{eq6} implies that all the inequalities \eqref{eq12}, \eqref{ret1} and \eqref{ret2} are equalities. Hence,
\begin{enumerate}
        \item $\xi_\eps^i = P_i\xi_\eps^i$ for all $\eps\neq (1,1,\ldots ,1)$ and all $1\leq i\leq m$ (cf. \eqref{eq22} and \eqref{eq21}), and
        \item by the Cauchy-Schwarz argument, for all $\eps\in\{0,1\}^m$ with more than one coordinate
        which is zero, $\|\xi_\eps^i\|\overset{({\rm i})}=\|P_i\xi_\eps^i\|=0$ for all $i$.
\end{enumerate}
Thus, if equality holds in \eqref{eq6}, then $\xi_i\in P_i(\CH)$ and
$\xi_i\in \bigoplus_{j\neq i} \CH_j$, as claimed.
$\endproof$

\vspace{.2cm}

\noindent {\it Proof of Proposition~\ref{thm1}}. \eqref{eq1} is the content of \cite[Lemma~2]{S}. Now, assume that equality holds in \eqref{eq1} and that $\Phi(\sum_{i\neq j}a_ix_i)$ is finite for all $j$. We are going to "backtrack" the proof of \cite[Lemma~2]{S} to show that \eqref{eq3} holds. We will assume that $\|x_j\|_2=1$ for all $j$.

With 
\[
\xi_j = b_j\CJ\Bigg({\textstyle \frac{1}{\sqrt{1-a_j^2}}}\sum_{i\neq j}a_ix_i\Bigg), \qquad 1\leq j\leq n+1,
\]
equality in \eqref{eq1} implies (cf. \cite[proof~of~Lemma~2]{S}) that 
\begin{equation}\label{eq4}
\Phi\Big(\sum_{j=1}^{n+1}a_jx_j\Big) = \Big\|\sum_{j=1}^{n+1}\xi_j\Big\|_2^2 = n\sum_{j=1}^{n+1}\|\xi_j\|_2^2.
\end{equation}
Let $M=W\cc(x_1,\ldots,x_{n+1})$. We now apply Lemma~\ref{commuting projections} to the projections $E_1,\ldots,E_{n+1}\in B(L^2(M))$ introduced in \cite[proof of Lemma~2]{S}. That is, $E_j$ is the projection onto $L^2(W\cc(x_1,\ldots,\hat{x_j},\ldots,x_{n+1}))$. Note that the subspace $\CH_j$ defined in Lemma~\ref{commuting projections}, 
\[
\CH_j= \{\xi\in L^2(M)\,|\,E_k\xi =\xi, \,k\neq j\,, E_j\xi =0\},
\]
is in this case exactly $L^2(W\cc(x_j))$. Thus, the second identity in \eqref{eq4} and the fact that $\xi_j\bot \C\unit$, implies that 
\begin{equation}\label{eq5}
\xi_j \in \bigoplus_{i\neq j} (L^2(W\cc(x_i))\ominus \C\unit).
\end{equation}
With $E:L^2(M)\rightarrow L^2(M)$ the projection onto $L^2(W\cc(\sum_j{a_jx_j}))$ we have (cf. \cite[proof of Lemma~2]{S}):
\begin{equation}\label{eq7}
\CJ\Big(\sum_{j=1}^{n+1}a_jx_j\Big)= E\Big(\sum_{j=1}^{n+1}\xi_j\Big).
\end{equation}
The first identity in \eqref{eq4} then implies that $E\Big(\sum_{j=1}^{n+1}\xi_j\Big)=\sum_{j=1}^{n+1}\xi_j$, and so
\[
\CJ\Big(\sum_{j=1}^{n+1}a_jx_j\Big) = \sum_{j=1}^{n+1}\xi_j \in \bigoplus_{i=1}^{n+1} (L^2(W\cc(x_i))\ominus \C\unit).
\]
Now choose elements $\eta_j\in L^2(W\cc(x_j))\ominus \C\unit$, $1\leq j\leq n+1$, such that
\begin{equation}\label{eq20}
\CJ\Big(\sum_{j=1}^{n+1}a_jx_j\Big) = \sum_{j=1}^{n+1} \eta_j.
\end{equation}
Then
\[
0= \Big[\sum_{i=1}^{n+1}a_ix_i, \sum_{j=1}^{n+1} \eta_j\Big] 
= \sum_{i\neq j} \big(a_ix_i \eta_j - \eta_i a_jx_j\big).
\]
A standard application of freeness shows that for $(i,j)\neq (k,l)$, the terms $a_ix_i \eta_j - \eta_i a_jx_j$ and $a_kx_k \eta_l - \eta_k a_lx_l$ are perpendicular elements of $L^2(M)$. Thus, the above identity implies that for all $i\neq j$,
\begin{equation}\label{eq8}
a_ix_i \eta_j = a_j\eta_i x_j.
\end{equation}

With $L^2(W\cc(x_j))^0 = L^2(W\cc(x_j))\ominus \C\unit$, $1\leq j\leq n+1$, consider the free product of Hilbert spaces
\[
\C\unit \oplus\Bigg(\bigoplus_{p\geq 1}\Big(\bigoplus_{1\leq i_1,\ldots, i_p\leq n+1,\; i_1\neq i_2\neq\cdots\neq i_p} L^2(W\cc(x_{i_1}))^0\otimes L^2(W\cc(x_{i_2}))^0\otimes\cdots\otimes L^2(W\cc(x_{i_p}))^0\Big)\Bigg),
\]
and notice that $x_i\in L^2(W\cc(x_i))^0$ and $\eta_j\in L^2(W\cc(x_j))^0$.
It follows from unique decomposition within the free product that there is only one way that \eqref{eq8} can be fulfilled, namely when $\eta_j$ is proportional to $x_j$. That is, there exist $c_1,\ldots, c_{n+1}\in\R$ such that $\eta_j=c_jx_j$ and hence,
\begin{equation}\label{eq9}
\CJ\Big(\sum_{j=1}^{n+1}a_jx_j\Big) = \sum_{j=1}^{n+1} c_jx_j.
\end{equation}
We can assume that $a_1,\ldots,a_{n+1}>0$, and then by \eqref{eq8},
\[
c_j = \frac{c_1a_j}{a_1}, \qquad 1\leq j\leq n+1.
\]
In particular, all the $c_j$'s have the same sign. Taking inner product with $\sum_{j=1}^{n+1}a_jx_j$ in \eqref{eq9}, we find that
\begin{equation}\label{eq10}
\sum_{j=1}^{n+1} a_jc_j =1,
\end{equation}
so that the $c_j$'s must be positive. Also, since $\sum_ja_j^2=1$, we have that $\sum_j c_j^2 \geq 1$. But
\[
\sum_{j=1}^{n+1} c_j^2 = \frac{c_1^2}{a_1^2},
\]
and so $c_1\geq a_1$, and in general, $c_j\geq a_j$. Then by \eqref{eq10}, $c_j=a_j$, and \eqref{eq3} holds. As mentioned in the introduction, this implies that $\sum_{j=1}^{n+1}a_jx_j$ is (0,1)-semicircular (when $\|x_1\|_2=1$). $\endproof$

\vspace{.2cm}

\begin{cor}\label{cor1} Let $x_1,\ldots, x_{n+1}$ be as in Proposition~\ref{thm1} and let $a_1,\ldots, a_{n+1}\in \R$ with $\sum_ja_j^2=1$. Then
\begin{equation}\label{eq11}
\chi\Big(\sum_{j=1}^{n+1}a_jx_j\Big) \geq \sum_{j=1}^{n+1} \frac{1-a_j^2}{n} \chi\Bigg({\textstyle \frac{1}{\sqrt{1-a_j^2}}}\sum_{i\neq j} a_ix_i\Bigg).
\end{equation}
Moreover, if $\chi(\sum_{i\neq j} a_ix_i) >-\infty$ for all $j$, then equality in \eqref{eq11} implies that $\sum_j a_jx_j$ is semicircular. 
\end{cor}

\noindent \proof The inequality \eqref{eq11} was proven by D.~Shlyakhtenko in \cite[Theorem~2]{S}. Now, assume that $\chi(\sum_{i\neq j} a_ix_i) >-\infty$ for all $j$ and that 
\[
\chi\Big(\sum_{j=1}^{n+1}a_jx_j\Big) = \sum_{j=1}^{n+1} \frac{1-a_j^2}{n} \chi\Bigg({\textstyle \frac{1}{\sqrt{1-a_j^2}}}\sum_{i\neq j} a_ix_i\Bigg).
\]
Take (0,1)-semicirculars $s_1,\ldots,s_{n+1}$ such that $x_1,\ldots, x_{n+1}, s_1,\ldots,s_{n+1}$ are free, and put 
\[
x_j^{(t)}= x_j + \sqrt{t}\, s_j.
\]
Then by assumption,
\begin{equation}\label{eq13}
\int_0^\infty \Bigg[\sum_{j=1}^{n+1}\frac{1-a_j^2}{n} \Phi\Bigg({\textstyle \frac{1}{\sqrt{1-a_j^2}}}\sum_{i\neq j} a_ix_i^{(t)}\Bigg)- \Phi\Big(\sum_{j=1}^{n+1}a_jx_j^{(t)}\Big)\Bigg] \d t =0.
\end{equation}
Applying Proposition~\ref{thm1} with $b_j=\frac1n \sqrt{1-a_j^2}$, we see that the integrand in 
\eqref{eq13} is positive. Thus, \eqref{eq13} can only be fulfilled if for a.e. $t> 0$,
\begin{equation}\label{eq14}
\sum_{j=1}^{n+1}\frac{1-a_j^2}{n} \Phi\Bigg({\textstyle \frac{1}{\sqrt{1-a_j^2}}}\sum_{i\neq j} a_ix_i^{(t)}\Bigg)= \Phi\Big(\sum_{j=1}^{n+1}a_jx_j^{(t)}\Big).
\end{equation}
In fact, since both sides of \eqref{eq14} are right continuous
functions of $t$ (cf. \cite{V2}), we have equality for all $t>
0$. Then by Proposition~\ref{thm1}, $\sum_{j=1}^{n+1}a_jx_j^{(t)}$ is
semicircular. By additivity of the $\CR$-transform, this can only
happen if $\sum_{j=1}^{n+1}a_jx_j$ is semicircular. $\endproof$ 

\vspace{.2cm}

\noindent {\it Proof of Theorem~\ref{thm2}.} The inequality \eqref{eq2} was proven by D.~Shlyakhtenko in \cite{S}. Now, assume that $\chi(x_1)>-\infty$ and that 
\[
\chi\Big(\frac{x_1+\cdots+x_{n+1}}{\sqrt{n+1}}\Big) = \chi\Big(\frac{x_1+\cdots+x_n}{\sqrt n}\Big).
\]
If we replace $x_j$ by $\frac{x_j-\tau(x_j)}{\|x_j-\tau(x_j)\|_2}$, we will still have equality. Hence, we will assume that $\tau(x_j)=0$ and that $\|x_j\|_2=1$. 
Now, 
\[
\chi\Big(\frac{x_1+\cdots+x_{n+1}}{\sqrt{n+1}}\Big) = \frac{1}{n+1}\sum_{j=1}^{n+1}\chi\Big(\frac{1}{\sqrt n}\sum_{i\neq j}x_i\Big),
\]
and by application of Corollary~\ref{cor1} with $a_j = \frac{1}{\sqrt{n+1}}$, $\frac{x_1+\cdots+x_{n+1}}{\sqrt{n+1}}$ must be semicircular. Additivity of the $\CR$-transform tells us that this can only happen if $x_1$ is semicircular. $\endproof$

\vspace{.2cm}

\noindent We would like to thank Serban Belinschi for pointing out to
us the following consequence of Theorem~\ref{thm2}:

\begin{cor}\label{stable} Among the freely stable compactly supported probability
  measures on $\R$, the semicirle laws are the only ones with finite
  free entropy.
\end{cor}

\noindent\proof By definition, a compactly supported probability
  measure $\mu$ on $\R$ is freely stable if for all $n\in\N$, there
  exist $a_n>0$, $b_n\in\R$, such that if $x_1,\ldots,x_n$ are freely
  independent self-adjoint elements which are distributed according to
  $\mu$, then
  \[
  \frac{1}{a_n}(x_1+\cdots+x_n)+ b_n
  \]
  has distribution $\mu$. Note that the set of freely stable laws is
  invariant under transformations by the affine maps
  $(\phi_{s,r})_{s\in\R, r>0}$, where
\[
\phi_{s,r}(t)=\frac{t-s}{r},\qquad t\in\R.
\]
Also, by \cite[p.~27]{VDN}, the semicirle laws are freely stable.

Suppose now that $\mu$ is a freely stable compactly supported
probability measure on $\R$. By the above remarks, we can assume that
$\mu$ has first moment 0 and variance 1. 

Let $x_1,x_2$ be freely independent
self-adjoint elements in  distributed according to
$\mu$. Since $\mu$ is freely stable, $\frac{x_1+x_2}{\sqrt 2}$ has
distribution $\mu$ as well (by the assumptions on $\mu$, $a_2=\sqrt2$
and $b_2=0$). But then
\[
\chi\Big(\frac{x_1+x_2}{\sqrt 2}\Big)= \chi(x_1),
\]
and by Theorem~\ref{thm2}, either $\chi(x_1)=-\infty$, or $x_1$ is semicircular. $\endproof$

\section{The Classical Case.}

In this section we are going to prove the classical analogue of Theorem~\ref{thm2}:

\begin{thm}\label{classical} Let $n\in\N$, and let $X_1,\ldots,X_{n+1}$ be i.i.d. random variables. Then
\begin{equation}\label{eqc1}
H\Bigg(\frac{X_1+\cdots+X_{n+1}}{\sqrt{n+1}}\Bigg) \geq  H\Bigg(\frac{X_1+\cdots+X_{n}}{\sqrt{n}}\Bigg).
\end{equation}
Moreover, if $H(X_1)>-\infty$ and if \eqref{eqc1} is an equality, then $X_1$ is Gaussian. 
\end{thm}
 
\vspace{.2cm} 

\begin{lemma}\label{Hermite} Let $n\in\N$. Then for every $m\in\N$, the $m$'th Hermite polynomial, $H_m$, satisfies:
\begin{equation}\label{eqc10}
n^{\frac m2}H_m\Big({\textstyle \frac{x_1+\cdots+x_n}{\sqrt n}}\Big)= \sum_{k_1,\ldots, k_n\geq 0,\, \sum_j k_j = m} \frac{m!}{k_1!k_2!\cdots k_n!}H_{k_1}(x_1)H_{k_2}(x_2)\cdots H_{k_n}(x_n).
\end{equation}
\end{lemma}

\noindent {\it Sketch of proof.} \eqref{eqc10} holds for $m=0$ ($H_0=1$) and for $m=1$ ($H_1(x)=2x$). Now, for general $m\in\N$,
\[
H_{m+1}(x)= 2xH_m(x)-2mH_{m-1}(x).
\]
\eqref{eqc10} then follows by induction over $m$. $\endproof$

\vspace{.2cm}

\begin{lemma}\label{L2-inclusion} Let $\mu\in{\rm Prob}(\R)$ be absolutely continuous w.r.t. Lebesgue measure, and let $\sigma_t\in{\rm Prob}(\R)$ denote the Gaussian distribution with mean 0 and variance $t$. Then if $\mu((-\infty,0])\neq 0$ and $\mu([0,\infty))\neq 0$, the following inclusion holds:
\begin{equation}\label{eqc5}
L^2(\R, \mu\ast\sigma_t)\subseteq L^2(\R,\sigma_t).
\end{equation}
\end{lemma}

\noindent \proof Let $f\in L^1(\R)$ denote the density of $\mu$ w.r.t. Lebesgue
measure. Then the density of $\mu\ast\sigma_t$ is given by
\begin{eqnarray*}
\frac{\d (\mu\ast\sigma_t)}{\d s}(s)&=& \frac{1}{\sqrt{2\pi t}}\Bigg(\int_{-\infty}^\infty f(u)\cdot e^{-\frac{u^2}{2t}}\cdot e^{\frac{su}{t}}\d u \Bigg)\cdot e^{-\frac{s^2}{2t}}\\
&=& \phi(s) \cdot \frac{\d \sigma_t}{\d s}(s),
\end{eqnarray*}
where
\begin{equation}
\phi(s) = \int_{-\infty}^\infty f(u)\cdot e^{-\frac{u^2}{2t}}\cdot e^{\frac{su}{t}}\d u.
\end{equation}
It follows that if $\phi$ is bounded away from 0, then \eqref{eqc5} holds. For $s\geq 0$ we have that
\begin{eqnarray*} 
\phi(s) &\geq & \int_0^\infty f(u)\cdot e^{-\frac{u^2}{2t}}\cdot e^{\frac{su}{t}}\d u\\
&\geq &  \int_0^\infty f(u)\cdot e^{-\frac{u^2}{2t}}\d u,
\end{eqnarray*}
and similarly for $s\leq0$:
\begin{eqnarray*} 
\phi(s) &\geq & \int_{-\infty}^0 f(u)\cdot e^{-\frac{u^2}{2t}}\d u.
\end{eqnarray*}
Since $\int_{-\infty}^0f(u)\d u >0$ and $\int_0^\infty f(u)\d u >0$, both of the integrals $\int_0^\infty f(u)\cdot e^{-\frac{u^2}{2t}}\d u$ and $\int_{-\infty}^0 f(u)\cdot e^{-\frac{u^2}{2t}}\d u$ are strictly positive. This completes the proof. $\endproof$

\vspace{.2cm}

\noindent {\it Proof of Theorem~\ref{classical}.} The inequality \eqref{eqc1}
was proven in \cite{ABBN}. Now, suppose $H(X_1)>-\infty$ and that
\eqref{eqc1} is an equality. We can assume that $X_1$ has first moment
$0$ and second moment 1. Take Gaussian random variables $G_1,\ldots,
G_{n+1}$ of mean 0 and variance 1 such that $X_1,\ldots, X_{n+1},
G_1,\ldots, G_n, G_{n+1}$ are independent. Then with
\[
X_j^{(t)}= X_j + \sqrt t\, G_j,
\]
\begin{equation}
H\Big({\textstyle \frac{X_1+\cdots+X_{n+1}}{\sqrt{n+1}}}\Big) =\frac12 \int_0^\infty \Big[\frac{1}{1+t}- \Big\|j\Big({\textstyle \frac{X_1^{(t)}+\cdots+X_{n+1}^{(t)}}{\sqrt{n+1}}}\Big)\Big\|_2^2\Big]\d t + \frac12 \log(2\pi e),
\end{equation}
where
\begin{equation}
j\Big({\textstyle \frac{X_1^{(t)}+\cdots+X_{n+1}^{(t)}}{\sqrt{n+1}}}\Big)= \Big(\frac{\d}{\d x}\Big)\cc (\unit) \in L^2\Bigg(\R, \mu_{\frac{X_1^{(t)}+\cdots+X_{n+1}^{(t)}}{\sqrt{n+1}}}\Bigg)
\end{equation}
is the score function. Since $X_1$ has mean 0 and finite entropy, $\mu_{X_1}$ and $\mu_{\frac{X_1+\cdots+X_{n+1}}{\sqrt{n+1}}}$ satisfy the conditions of Lemma~\ref{L2-inclusion}.

For $t>0$, define $f^{(t)}\in L^2(\R^{n+1},\otimes_{j=1}^{n+1}\mu_{X_j^{(t)}})$ by
\[
f^{(t)}(x_1,\ldots,x_{n+1})= j\Big({\textstyle \frac{X_1^{(t)}+\cdots+X_{n+1}^{(t)}}{\sqrt{n+1}}}\Big)\Big({\textstyle \frac{x_1+\cdots+x_{n+1}}{\sqrt{n+1}}}\Big).
\]
As in the free case (cf. \eqref{eq20}) equality in \eqref{eqc1} implies that for each $t>0$ there exists a function $g^{(t)}\in L^2(\mu_{X_1^{(t)}})$ such that $\int g^{(t)}\d \mu_{X_1^{(t)}}=0$ and 
\begin{equation}\label{eqc3}
f^{(t)}(x_1,\ldots,x_{n+1}) = \sum_{j=1}^{n+1} g^{(t)}(x_j).
\end{equation}
Because of Lemma~\ref{L2-inclusion} we can now write things in terms
of the Hermite polynomials $(H_m)_{m=0}^\infty$. That is, there exist
scalars $(\alpha_m)_{m=1}^\infty$ and $(\beta_m)_{m=1}^\infty$ such that
\[
f^{(1)}(x_1,\ldots,x_{n+1}) = \sum_{m=1}^\infty \alpha_m H_m\Big({\textstyle \frac{x_1+\cdots+x_{n+1}}{\sqrt{n+1}}}\Big),
\]
and
\[
g^{(1)}(x)=\sum_{m=1}^\infty \beta_m H_m(x).
\]
By Lemma~\ref{Hermite}, this implies that 
\begin{equation}\label{eqc4}\begin{split}
\sum_{j=1}^{n+1}&\sum_{m=1}^\infty \beta_m H_m(x_j) = \\
&\sum_{m=1}^\infty \frac{\alpha_m}{(n+1)^{\frac m2}}\sum_{\begin{array}{c} {\textstyle k_1,\ldots, k_{n+1}\geq 0},\\ {\textstyle \sum_j k_j = m}\end{array}} \frac{m!}{k_1!k_2!\cdots k_{n+1}!}H_{k_1}(x_1)H_{k_2}(x_2)\cdots H_{k_{n+1}}(x_{n+1}).\end{split}
\end{equation}
The functions $(H_{k_1}(x_1)H_{k_2}(x_2)\cdots
H_{k_{n+1}}(x_{n+1}))_{k_1,\ldots,k_{n+1}\geq 0}$ are mutually perpendicular in $L^2(R^{n+1},\otimes_{j=1}^{n+1} \sigma_1)$. Fix $m\geq 2$, and take $k_1,\ldots, k_{n+1}$ with $\sum_j k_j=m$ and $k_j \geq 1$ for at least two $j$'s. Then take inner product with $H_{k_1}(x_1)H_{k_2}(x_2)\cdots H_{k_{n+1}}(x_{n+1})$ on both sides of \eqref{eqc4} to see that $\alpha_m$ must be zero. That is, 
\[
j\Big({\textstyle \frac{X_1^{(1)}+\cdots+X_{n+1}^{(1)}}{\sqrt{n+1}}}\Big)\Big({\textstyle \frac{x_1+\cdots+x_{n+1}}{\sqrt{n+1}}}\Big) = \alpha_1 H_1\Big({\textstyle \frac{x_1+\cdots+x_{n+1}}{\sqrt{n+1}}}\Big)= 2\alpha_1{\textstyle \frac{x_1+\cdots+x_{n+1}}{\sqrt{n+1}}}.
\]
Since the score function of a random variable $X$, $j(X)$, satisfies
$\<j(X),X\>_{L^2(\mu_X)}=1$, we have that $\alpha_1=\frac12$, and so
\[
j\Big({\textstyle \frac{X_1^{(1)}+\cdots+X_{n+1}^{(1)}}{\sqrt{n+1}}}\Big)\Big({\textstyle \frac{x_1+\cdots+x_{n+1}}{\sqrt{n+1}}}\Big)= {\textstyle \frac{x_1+\cdots+x_{n+1}}{\sqrt{n+1}}}.
\]
Then ${\textstyle \frac{X_1^{(1)}+\cdots+X_{n+1}^{(1)}}{\sqrt{n+1}}}$
has Fisher information 1, implying that it is standard Gaussian. As in
the free case, using additivity of the logarithm of the Fourier
transform, this can only happen if $X_1$ is Gaussian. $\endproof$

{\small

\vspace{.2cm}

\noindent Hanne~Schultz,\\
Department of Mathematics and Computer Science,\\
University of Southern Denmark,\\
Campusvej 55,\\
5230 Odense M, Denmark\\
{\texttt schultz@imada.sdu.dk}

\end{document}